\documentclass[12pt,reqno]{article}
\usepackage{hyperref}

\let\eps=\varepsilon

\usepackage{amsmath}
\usepackage{amssymb}
\usepackage{amsthm}
\usepackage{srcltx}
\usepackage{enumerate}

\newcommand\N{{\mathbb{N}}}
\newcommand\R{{\mathbb{R}}}

\newcommand{\ds}{\displaystyle}
\newcommand{\txs}{\textstyle}

\newcommand{\be}{\begin{equation}}
\newcommand{\ee}{\end{equation}}
\newcommand{\ba}{\begin{array}}
\newcommand{\ea}{\end{array}}
\newcommand{\bea}{\begin{eqnarray}}
\newcommand{\eea}{\end{eqnarray}}
\newcommand{\Bea}{\begin{eqnarray*}}
\newcommand{\Eea}{\end{eqnarray*}}
\newcommand{\bt}{\begin{Theorem}}
\newcommand{\et}{\end{Theorem}}

\newcommand{\bpr}{\begin{Proposition}}
\newcommand{\epr}{\end{Proposition}}
\newcommand{\bpb}{\begin{Problem}}
\newcommand{\epb}{\end{Problem}}

\newcommand{\bl}{\begin{Lemma}}
\newcommand{\el}{\end{Lemma}}
\newcommand{\bi}{\begin{itemize}}
\newcommand{\ei}{\end{itemize}}
\numberwithin{equation}{section}

\newtheorem{lemma}{Lemma}[section]
\newtheorem{theorem}[lemma]{Theorem}
\newtheorem{definition}[lemma]{Definition}
\newtheorem{proposition}[lemma]{Proposition}
\newtheorem{problem}[lemma]{Problem}

\theoremstyle{definition}

\def\2to{\to\!\!\!\!\!\!{}_{_2}\;\;}
\def\w2to{\rightharpoonup\!\!\!\!\!\!{}_{_2}\;\;}
\def\ws2to{\rightharpoonup\!\!\!\!\!\!{}_{_2}\!\!\!{}^{*}\;\;\;}
\def\wto{\rightharpoonup}
\def\wsto{\rightharpoonup\!\!\!\!\!\!{}^{*}\;\,}

\title{Variational Approach to Homogenization of
Doubly-Nonlinear Flow in a Periodic Structure}

\author{A. K. Nandakumaran
\thanks{Department of Mathematics, Indian Institute of Science,
Bangalore 560012, India --
\hfill\break
email: nands@math.iisc.ernet.in } ,
Augusto Visintin
\thanks{Dipartimento di Matematica dell'Universit\`a degli Studi di Trento
via Sommarive 14,  38050 Povo di Trento, Italia -- email: Visintin@science.unitn.it }
}

\date{October 6, 2014}

\begin{document}
\maketitle 

\begin{abstract}
This work deals with the homogenization of an initial- and boundary-value problem
for the doubly-nonlinear system
\begin{eqnarray}
&\ds D_t w -\nabla\cdot \vec z = \nabla\cdot \vec h(x,t,x/\eps)
\label{sis1=}
\\ [1mm]
&\ds w\in \alpha(u,x/\eps)
\label{sis2=}
\\ [1mm]
&\ds \vec z\in \vec\gamma(\nabla u,x/\eps).
\label{sis3=}
\end{eqnarray}
Here $\eps$ is a positive parameter, and the prescribed mappings
$\alpha$ and $\vec\gamma$ are maximal monotone
with respect to the first variable and periodic with respect to the second one.

The inclusions \eqref{sis2=} and \eqref{sis3=} are here formulated as
{\it null-minimization principles,\/} via the theory of Fitzpatrick [MR 1009594].
As $\eps\to 0$, a two-scale formulation is derived via Nguetseng's notion
of two-scale convergence, and a (single-scale) homogenized problem is then retrieved.
\end{abstract}

\bigskip
\noindent{\bf AMS Classification (2000): }
35B27, 
35K60, 
49J40, 
78M40. 


\bigskip
\section{Introduction}

This paper deals with the homogenization of a class of
doubly-nonlinear parabolic equations of the form
\be \label{pb}
\left\{
\ba{ll}
D_t w_\eps -\nabla\cdot \vec z_\eps = \nabla\cdot \vec h(x,t,x/\eps)
\\ [2mm]
w_\eps\in \alpha(u_\eps,x/\eps)
\\ [2mm]
\vec z_\eps\in \vec\gamma(\nabla u_\eps,x/\eps)
\ea\right.
\qquad \mbox{ in }\Omega \times {}]0,T[.
\ee
Here $\Omega$ is a bounded domain of $\R^N$, $T>0$,
and $\eps$ is a positive parameter. The mappings
\be
\alpha:\R \times \R^N \to \mathcal{P}(\R),
\qquad
\vec\gamma: \R^N\times \R^N\to \mathcal{P}(\R^N)
\ee
are prescribed, and are maximal monotone with respect to the first argument
and periodic with respect to the second one.
The known source field $\vec h$ is also periodic with respect to the third argument.
We also assume that
\be\label{bc}
\ba{ll}
u_\eps = 0 \qquad \mbox{ on }\partial\Omega \times {}]0,T[,
\\ [2mm]
w_\eps(x,0) = w^0(x,x/\eps) \qquad \mbox{ for }x\in\Omega,
\ea
\ee
for a prescribed periodic function $w^0$.
All periods are assumed to coincide.

Problems of the form \eqref{pb} arise in several physical contexts:
e.g., this may represent the entropy balance in diffusion phenomena;
$\alpha$ may be the subdifferential of a dissipation potential.
Existence of a solution for an associated boundary- and initial-value problem
was proved e.g.\ by DiBenedetto and Showalter \cite{DiSh} 
and by Alt and Luckhaus \cite{AlLu}. 

In the case of single-valued operators, the homogenization of a
similar system was already studied by H. Jian \cite{Ji}.
This was also used to model filtration in porous media by
A.K.N. and M. Rajesh \cite{NaRa1}, \cite{NaRa2}, \cite{NaRa3}.
More precisely, in \cite{NaRa1} a quasi-linear equation of the form
$$
\partial_{t} \alpha(u_\varepsilon,x/\eps) - \nabla\cdot
\vec\gamma(u_\varepsilon,\nabla u_\varepsilon, x/\eps, t/\eps) = h(x,t)
$$
was considered with appropriate boundary and initial conditions, thus also
accounting for high-frequency oscillations with respect to time. 
The same equation was also addressed by A.K.N. and M. Rajesh
\cite{NaRa2}, \cite{NaRa3}, dealing with in a porous medium with
Neumann and Dirichlet boundary conditions, respectively. In
\cite{NaRa1}, \cite{NaRa2} two-scale convergence was used extensively. 
It should be noticed that the Dirichlet condition 
on the boundary of the holes may yield different homogenized problems, 
that depend on the asymptotic relation between 
the size of the holes and the period $\varepsilon$.

The homogenization of quasi-linear equations has been studied by various authors, see 
e.g.\ \cite{AvLi}, \cite{BrFrMu}, \cite{FuMo}, \cite{OlKoZh}. 
The homogenization of doubly-nonlinear equations of the form \eqref{pb} occurring in electromagnetic processes in composites and in Stefan-type problems 
was performed in \cite{ViSIMA07}, \cite{ViCom08}.

\medskip
Each of the inclusions \eqref{pb}$_2$ and \eqref{pb}$_3$
is equivalent to a variational inequality.
On the basis of the Fitzpatrick theory \cite{Fi},
here we convert the system above into the coupling of a linear PDE with a
{\it null-minimization\/} problem, along the lines of \cite{ViDNE}.
We then study the limit behaviour for vanishing $\eps$.

This note is organized as follows.
First in section 2 we briefly outline the Fitzpatrick theory
for the variational representation of maximal monotone operators.
In section 3 we describe the homogenization problem to be studied,
and in section 4 we prove existence of a solution.
We then let $\eps$ vanish.
In section 5 we derive the two-scale problem,
and in section 6 we then retrieve a single-scale system,
by this proving the desired homogenization theorem.
Finally in an appendix we briefly review Nguetseng's theory of two-scale convergence
and related properties of integral functionals;
these also include a result of \cite{ViMon} on
the homogenization of maximal monotone operators.

The novelty of this work stays in the use of a Fitzpatrick-type formulation
for homogenization, and in the derivation of a two-scale problem
as an intermediate steps towards homogenization. 

The results of this note may be extended in several directions;
for instance explicit dependence on time may be assumed in the nonlinear operator,
and time-homogenization may also be considered.
The homogenization of several other quasilinear equations may also be studied,
including doubly-nonlinear systems of the form
\begin{eqnarray}
&\ds w_\eps -\nabla\cdot \vec z_\eps = \nabla\cdot \vec h(x,t,x/\eps)
\label{sis1'}
\\ [1mm]
&\ds w_\eps\in \alpha(D_t u_\eps,x/\eps)
\label{sis2'}
\\ [1mm]
&\ds \vec z_\eps\in \vec\gamma(\nabla u_\eps,x/\eps),
\label{sis3'}
\end{eqnarray}
with $\alpha$ and $\vec\gamma$ as above.
Existence of a solution for an associated boundary- and initial-value problem 
was proved in \cite{CoVi}.

\bigskip
\section{Preliminaries} 

\noindent
In this section we illustrate the tenets of the Fitzpatrick theory
on the variational representation of maximal monotone operators,
that is at the basis of the procedures of the present work.
We also illustrate an idea of Brezis, Ekeland and Nayroles
for the variational formulation of monotone flows.
We refer e.g.\ to \cite{ViCalVar13} for a more detailed review.

\subsection{Variational representation of maximal monotone operators}

Let us first recall the Fenchel system,
which is a basic result of the theory of convex analysis,
see e.g.\ \cite{EkTe}, \cite{Ro}.
Let $V$ be a separable and reflexive real Banach space with dual $V'$, let
$\psi: V\to\R\cup \{+\infty\}$
be a convex and lower semicontinuous function, and
$\psi^*: V'\to\R \cup \{+\infty\}$ be its conjugate function, namely,
\be
\psi^*(v'):= \sup_{v\in V} \;\{\langle v',v\rangle - \psi(v)\}
\qquad \forall v'\in V'.
\ee
It is known that $\psi$, $\psi^*$ and the subdifferential $\partial \psi$
satisfy the following Fenchel system:
\be\label{fenchel}
\left\{
\ba{ll} \psi(v) + \psi^*(v')\geq \langle v',v\rangle
\qquad\forall (v,v') \in V\times V',
\\ [2mm]
\psi(v) + \psi^*(v') =  \langle v',v\rangle
\mbox{ \ \ if and only if \ \ }
v' \in \partial \psi(v).
\ea\right.
\ee

Let now $\alpha: V\times {\cal P}(V')$ be a multivalued mapping.
In \cite{Fi} Fitzpatrick introduced the following convex and lower semicontinuous function:
\be\label{fi}
\ba{lll}
f_\alpha(v,v')
&:= \langle v',v\rangle +\sup \;\{ \langle v'-v_0',v_0-
v\rangle: \forall v_0' \in \alpha(v_0)\}
\\ [2mm]
&= \sup \;\{\langle v',v_0\rangle - \langle v_0',v_0- v\rangle:
\forall v_0' \in \alpha(v_0)\}
\ea\ee
for all $(v,v') \in V\times V'$, and proved that, whenever $\alpha$ is maximal monotone,
\be\label{fitz1}
\left\{\ba{ll} f_\alpha(v,v')\geq \langle v',v\rangle
\qquad\forall (v,v') \in V\times V',
\\ [2mm]
f_\alpha(v,v')= \langle v',v\rangle
\mbox{ \ \ if and only if \ \ } v' \in \alpha(v).
\ea\right.
\ee
This system obviously extends \eqref{fenchel}.
Nowadays $f_\alpha$ is called the {\it Fitzpatrick function\/} of $\alpha$.

The inclusion $v' \in \alpha(v)$ is thus equivalent to
\be\label{nullmin}
f_\alpha(v,v') - \langle v',v\rangle
= \inf\; \{f_\alpha(r,r') - \langle r',r\rangle: (r,r') \in V\times V'\} =0,
\ee
that we label as a {\it null-minimization problem.\/}

Next we review the notion of {\it (variational) representation\/} of monotone operators.

\begin{definition}
We shall say that a lower semicontinuous convex function
$f: V\times V' \to \R \cup \{+\infty\}$ (variationally) represents a 
(necessarily monotone) operator $\alpha: V\to\mathcal{P}(V')$ in the sense of Fitzpatrick, whenever
\be\label{fitz2}
\left\{\ba{ll}
f(v,v')\geq \langle
v',v\rangle \qquad \forall (v,v') \in V\times V',
\\ [2mm]
f(v,v')= \langle v',v\rangle
\mbox{ \ \ if and only if \ \ } v' \in \alpha(v).
\ea\right.\ee
\end{definition}

Such a function is called a {\it representative function.\/}
For instance, because of (\ref{fi})--(\ref{nullmin}),
$\alpha$ is represented by the function $f_\alpha$.
If $\alpha= \partial \psi$, then because of \eqref{fenchel} $\alpha$ is also represented
by the {\it Fenchel function\/} $g_\alpha(v,v') := \psi(v) + \psi^*(v')$.

\subsection{The Brezis-Ekeland-Nayroles variational formulation of flows}

Let us assume that we are given a triplet of (real) Banach spaces
\be
V\subset H=H'\subset V' \mbox{ \ with continuous and dense injections.}
\ee
On the basis of the Fenchel system \eqref{fenchel}, under suitable restrictions,
for any prescribed lower semicontinuous and convex function
$\psi:V\to \R \cup \{+\infty\}$, any $u^*\in L^2(]0,T[;V')$ and any $u^0\in H$,
Brezis and Ekeland \cite{BrEk} and Nayroles \cite{Nay} independently reformulated
the gradient flow
\be\label{flow1}
D_t u + \partial \psi(u)= u^* \qquad\mbox{ in }{}]0,T[
\ee
as the null-minimization of the functional
\be\label{BEN1}
\Phi_1(v, u^*) =
\int_0^T [\psi(v) + \psi(u^*- D_tv)] \, dt + {\txs{1\over2}} (\|v(T)\|_H^2 - \|u(0)\|_H^2)
- \langle u^*,v\rangle,
\ee
as $v$ ranges in $H^1(0,T;V') \cap L^2(0,T;V)$ ($\subset C^0([0,T];H)$).
More generally, see \cite{ViAMSA}, for any maximal monotone
$\alpha: V\to\mathcal{P}(V')$,
denoting by $f_\alpha$ a representative functions of $\alpha$, the monotone flow
\be\label{flow2}
D_t u + \alpha(u)= u^* \qquad\mbox{ in }{}]0,T[
\ee
may be represented as the null-minimization of the functional
\be\label{BEN2}
\Phi_2(v, u^*) =
\int_0^t f_\alpha(v, u^*- D_tv) \, dt + {\txs{1\over2}} (\|v(T)\|_H^2 - \|u(0)\|_H^2)
- \langle u^*,v\rangle.
\ee

\bigskip
\section{Weak Formulation of the $\eps$-Problem}
\setcounter{equation}{0}

In this section we provide two equivalent formulations of the system \eqref{pb}
coupled with appropriate initial- and boundary-conditions in a periodic medium.

Let $Y= {} ]0,1[^N$ be the unit cell, and let us assume that
\be\ba{ll} \label{defg}
&\ds g: \R\times Y\to \R \cup \{+\infty\}
\mbox{ \ is measurable w.r.t.\ }\mathcal{B}(\R) \otimes \mathcal{L}(Y),
\\ [2mm]
&\ds \varphi(\cdot,y) \mbox{ is convex and lower semicontinuous for a.e.\ $y$, }
\ea\ee
\be\label{g-growth}
\exists c_1, c_2 >0: \forall v\in \R, \qquad
|\varphi(v,y)|\leq c_1|v|^2 +c_2 \quad\mbox{ \ for a.e.\ }y\in Y.
\ee
By definition of the convex conjugate function $\varphi^*(\cdot,y)$, it follows that
\be\label{gstar-growth}
\exists L,M>0: \forall v\in \R, \qquad
|\varphi^*(v,y)| \geq L |v|^2 -M \quad\mbox{ \ for a.e.\ }y\in Y.
\ee
Let us set
\be\label{subdif}
\alpha(\cdot,y) = \partial\varphi(\cdot,y)
\qquad\mbox{ \ for a.e.\ }y\in Y;
\ee
by this we denote the {\it subdifferential\/}
with respect to the first variable (see e.g.\ \cite{EkTe}, \cite{Ro}).
The multivalued map $\alpha:\R\times Y \to {\cal P}(\R^N)$
is then measurable with respect to the $\sigma$-algebra
$\mathcal{B}(\R)\bigotimes\mathcal{L}(Y)$,
and $\alpha(\cdot, y)$ is maximal monotone for a.e.\ $y$.
Moreover $\alpha(v,\cdot)$ is measurable
for any measurable function $v: \R^N\times Y\to \R$. (See the Appendix.)

Let us assume that
\be\ba{ll}
&\ds \vec\gamma: \R^N \!\times\! Y \to \mathcal{P}(\R^N)
\mbox{ \ is measurable w.r.t. $\mathcal{B}(\R^N) \otimes \mathcal{L}(Y)$, }
\\ [2mm]
&\ds \vec\gamma(\cdot,y) \mbox{ \ is maximal monotone for a.e.\ }y,
\ea\ee
and that there exist nonnegative constants $k, a, b$ such that
\begin{eqnarray}
&|\vec z| \leq k(1+|\vec \zeta|)
\qquad \forall (\vec \zeta,\vec z)\in \mbox{graph} (\vec\gamma(\cdot,y)),
\mbox{ for a.e.\ }y,
\label{growth}
\\[2mm]
&\vec z\cdot \vec \zeta \geq a(| \vec z |^2 + |\vec \zeta|^2) - b
\qquad \forall (\vec \zeta , \vec z)\in \mbox{graph} (\vec\gamma(\cdot,y)),
\mbox{ for a.e.\ }y.
\label{bound-below}
\end{eqnarray}

Let us also assume that $\Omega$ is a bounded domain of $\R^N$ of Lipschitz class,
and that, setting $\Omega_T:= \Omega \times {}]0,T[$,
\be\ba{ll} \label{defh}
&\ds \vec h: \Omega_T\times Y\to \R^N
\mbox{ \ is measurable w.r.t.\ }\mathcal{B}(\Omega_T) \otimes \mathcal{L}(Y),
\\ [2mm]
&\ds \vec h(\cdot,\cdot,y) \in L^2(\Omega_T)^N \qquad\mbox{ for a.e.\ }y,
\ea\ee
\be\ba{ll} \label{defw}
&\ds w^0: \Omega\times Y\to \R
\mbox{ \ is measurable w.r.t.\ }\mathcal{B}(\Omega) \otimes \mathcal{L}(Y),
\\ [2mm]
&\ds w^0(\cdot,y) \in L^2(\Omega) \qquad\mbox{ for a.e.\ }y.
\ea\ee

We extend all of these functions $Y$-periodically to $\R^N$
with respect to the argument $y$, and set
\begin{eqnarray}
&\varphi_\eps (v,x):= \varphi(v,x/\eps)
\qquad
\forall v\in\R, \mbox{ for a.e.\ }x\in \R^N,
\label{index1}
\\ 
&\vec\gamma_\eps(v,x) := \vec\gamma(v,x/\eps)
\qquad
\forall v\in\R^N, \mbox{ for a.e.\ }x\in \R^N,
\label{index3}
\\
&\vec h_\eps(x,t) := \vec h(x,t,x/\eps)
\qquad
\mbox{ for a.e.\ }(x,t)\in \Omega_T,
\label{index4}
\\
&w^0_\eps(x) := w^0(x,x/\eps)
\qquad
\mbox{ for a.e.\ }x\in \Omega.
\label{index5}
\end{eqnarray} 
We shall deal with the homogenization of the following doubly-nonlinear system
\begin{eqnarray}
&D_tw_\eps -\nabla\cdot \vec z_\eps = \nabla\cdot \vec h_\eps
\qquad\mbox{ in ${\cal D}'(\Omega)$, a.e.\ in }{}]0,T[,
\label{inhomoeqn1a}
\\
&w_\eps\in \partial\varphi_\eps(u_\eps,x) \qquad\mbox{ a.e.\ in }\Omega_T,
\label{inhomoeqn1b}
\\
&\vec z_\eps\in \vec\gamma_\eps(\nabla u_\eps,x) \qquad\mbox{ a.e.\ in }\Omega_T,
\label{inhomoeqn1c}
\\
&u_\eps = 0 \qquad\mbox{ a.e.\ on }\partial\Omega \times {}]0,T[,
\label{inhomoeqn1d}
\\
&w_\eps(\cdot,0)= w^0_\eps \qquad\mbox{ a.e.\ in }\Omega.
\label{inhomoeqn1e}
\end{eqnarray}

We shall assume that 
\be\ba{ll}\label{bdd}
&\vec z_\eps \in L^2(0,T;V'),
\quad
w^0_\eps\in L^2(\Omega)
\\[2mm]
&\mbox{and are uniformly bounded w.r.t.\ $\eps$ in these spaces.}
\ea
\ee

Next we introduce the Hilbert triplet
\be
V= H^1_0(\Omega) \subset H = L^2(\Omega) = H' \subset V'= H^{-1}(\Omega)
\ee
(with continuous and dense injections),
and reformulate the system \eqref{inhomoeqn1a}--\eqref{inhomoeqn1e}
in weak form as follows, for any $\eps>0$.

\medskip
\begin{problem}\label{problem1}
Find $(u_\eps, w_\eps, \vec z_\eps) \in
L^2(0,T; V) \times L^2(\Omega_T)\times L^2(\Omega_T)^N$ such that
\begin{eqnarray}
&\ba{rr}
\ds \int\!\!\!\int_{\Omega_T} [(w^0_\eps - w_\eps) D_tv
+ (\vec z_\eps + \vec h_\eps) \cdot \nabla v] \, dxdt = 0
\\ [3mm]
\forall v\in H^1(0,T;V), v(\cdot, T)=0,
\ea
\label{inhomoeqn2a}
\\
&w_\eps \in \partial\varphi_\eps(u_\eps,x)
\qquad\mbox{ a.e.\ in }\Omega_T,
\label{inhomoeqn2b}
\\
&\vec z_\eps\in \vec\gamma_\eps(\nabla u_\eps,x)
\qquad\mbox{ a.e.\ in }\Omega_T.
\label{inhomoeqn2c}
\end{eqnarray}
\end{problem}

The equation \eqref{inhomoeqn2a} yields
\be \label{PDE=}
D_tw_\eps - \nabla\cdot \vec z_\eps = \nabla\cdot \vec h_\eps
\qquad\mbox{ in $V'$, a.e.\ in }]0,T[.
\ee
By comparing the terms of this equation we have
$D_tw_\eps \in  L^2(0,T; V')$, whence
\be
w_\eps \in H^1(0,T;V')\subset C^0([0,T];V')\quad
\mbox {(by an obvious identification). }
\ee
The equation \eqref{inhomoeqn2a} then also entails \eqref{inhomoeqn1e}.
Conversely, \eqref{inhomoeqn1e} and \eqref{PDE=} yield \eqref{inhomoeqn2a}.

\bigskip
Next we reformulate \eqref{inhomoeqn2b} and \eqref{inhomoeqn2c}
via the Fitzpatrick theorem \eqref{fitz1}.
First we denote by $f_{\vec\gamma_\eps}(\cdot, \cdot,x)$
a representative function of $\vec\gamma_\eps(\cdot,x)$ for a.e.\ $x$.
For any $(u,w,\vec z)$ that satisfies \eqref{inhomoeqn2a}, we set
\be \label{Phi-defn}\ba{ll}
&\ds\Phi_\eps (u,w, \vec z) :=
\\ [2mm]
&\ds \int\!\!\!\int_{\Omega_T} [\varphi_\eps(u,x) +
\varphi_\eps^*(w,x) -wu + f_{\vec\gamma_\eps}(\nabla u, \vec z,x)
-\nabla u\cdot \vec z ] \, dxdt, \ea\ee where $\varphi_\eps$ is
defined as in \eqref{index1}. We then define the
infinite-dimensional manifold \be X_\eps= \{(u, w, \vec z)\in
L^2(0,T; V) \times L^2(\Omega_T)\times L^2(\Omega_T)^N \mbox{ that
fulfill } \eqref{inhomoeqn2a}\}. \ee

For any $\eps>0$ we shall consider the following problem,
in which a PDE is coupled with a null-minimization problem.

\begin{problem}\label{problem2}
Find $(u_\eps, w_\eps, \vec z_\eps) \in X_\eps$ such that
\be\label{minimization1}
\Phi_\eps(u_\eps, w_\eps, \vec z_\eps) = \inf_{X_\eps} \Phi_\eps = 0.
\ee
\end{problem}

\begin{proposition}\label{equiv}
For any $\eps>0$, Problems~\ref{problem1} and \ref{problem2} are mutually equivalent.
\end{proposition}

\proof
By the Fenchel and Fitzpatrick systems \eqref{fenchel} and \eqref{fitz1},
the functional $\Phi_\eps$ is nonnegative.
The null-minimization \eqref{minimization1} is thus equivalent to the inequality
\be
\Phi_\eps(u_\eps, w_\eps, \vec z_\eps) \le 0,
\ee
or also to the system of the two inequalities
\begin{eqnarray}
&\ds \int\!\!\!\int_{\Omega_T}
[\varphi_\eps(u_ \eps,x) + \varphi_\eps^*(w_ \eps,x) -w_ \eps  u_ \eps]
\, dxdt\le 0,
\label{mini-inequality1}
\\
&\ds \int\!\!\!\int_{\Omega_T} [f_{\vec\gamma_\eps}(\nabla u_\eps , \vec z_\eps,x)
- \nabla u_\eps \cdot \vec z_\eps] \, dxdt \le 0.
\label{mini-inequality2}
\end{eqnarray}

The integrand of either functional is pointwise nonnegative,
so that by \eqref{fenchel} and \eqref{fitz1}
these inequalities are respectively equivalent to \eqref{inhomoeqn2b}
and \eqref{inhomoeqn2c}.
\qed

\bigskip
\section{Existence of a Solution of the $\eps$-Problems}

In this section we prove existence of a solution of Problem~\ref{problem1}.
Although this result has already been proved e.g.\ in \cite{DiSh},
here we present an argument based on the equivalence with Problem~\ref{problem2}.
We do so for the sake of completeness,
and also because this argument provides the uniform estimates that
we shall use in the homogenization procedure in the next section.

\subsection{Approximation by time-discretization}

Let us fix any $\eps >0$, any $m\in\N$, set $k=T/m$ and
\be
\vec h_{\eps m}^n =\ds\frac{1}{k}\int_{(n-1)k}^{nk} \vec h_\eps(\cdot,t) \, dt
\;(\in L^2(\Omega)^N) \qquad\mbox{ for }n = 1,...,m.
\ee
For any $\eps>0$ and any $m$, let us then consider the following
time-discretized problem.

\begin{problem}\label{problem3}
Find $(u^n_{\eps m}, w^n_{\eps m}, \vec z ^n_{\eps m}) \in V\times H\times H^N$
($n=1,...m$), such that, setting $w^0_{\eps m} = w^0_\eps$, for $n=1,...m$
\begin{eqnarray}
&w^n_{\eps m} -k\nabla\cdot \vec z ^n_{\eps m}
= w^{n-1}_{\eps m} + k \nabla\cdot \vec h_{\eps m}^n \qquad \mbox{ in } \, V',
\label{inhomoeqn-des1a}
\\
&w^n_{\eps m}\in \partial\varphi_\eps(u^n_{\eps m},x)
\qquad \mbox{ a.e.\ in }\Omega,
\label{inhomoeqn-des2a}
\\
&\vec z ^n_{\eps m}\in \vec\gamma_\eps(\nabla u^n_{\eps m},x)
\qquad \mbox{ a.e.\ in }\Omega.
\label{inhomoeqn-des1c}
\end{eqnarray}
\end{problem}

Defining
$\Lambda_{\eps,m}(v) = \partial\varphi_\eps(v,x)
- k\nabla\cdot \vec\gamma_\eps(\nabla v,x)
\in \mathcal{P}(V')$ for any $v\in V$,
the system \eqref{inhomoeqn-des1a}--\eqref{inhomoeqn-des1c} is equivalent to
\be\label{equiin}
\Lambda_{\eps,m} (u^n_{\eps m},x)\ni  w^{n-1}_{\eps m}
+  k\nabla \cdot \vec h_{\eps m}^n
\qquad \mbox{ in }V', n = 1,..., m.
\ee
By the assumptions \eqref{defg}--\eqref{defw}, for any $\eps,m$
the operator $\Lambda_{\eps,m}: V\to \mathcal{P}(V')$ is maximal monotone and coercive.
The inclusion \eqref{equiin} has then at least one solution, and this solves
Problem~\ref{problem3}.

Let us now define time-interpolate functions as follows.
For any family $\{v_m^n\} _{n=0,\dots,m}\subset \R$, let us denote
by $v_m$ the piecewise-linear time-interpolate of $v_m^0:=v^0$,
$v_m^1$, ..., $v_m^m$ a.e.\ in $\Omega$.
Let us denote by $\bar v_m$ the corresponding piecewise-constant interpolate
function, that is, $\bar v_m(t) :=v^n_m$ if $(n-1)k<t\le nk$ for $n=1,\dots,m$.

The system \eqref{inhomoeqn-des1a}--\eqref{inhomoeqn-des1c} then also reads
\begin{eqnarray}
&D_tw_{\eps m}-
\nabla\cdot \bar{\vec z}_{\eps m} = \nabla\cdot \bar{\vec h}_{\eps m}
\qquad\mbox{ \ in $V'$, a.e.\ in }{}]0,T[,
\label{inhomoeqn-des2a}
\\
&\bar{w}_\eps \in \partial\varphi_\eps(\bar{u}_{\eps m},x)
\qquad\mbox{ \ a.e.\ in }\Omega_T,
\label{inhomoeqn-des2b}
\\
&\bar{\vec z}_{\eps m}\in \vec\gamma_\eps(\nabla \bar{u}_{\eps m},x)
\qquad\mbox{ \ a.e.\ in }\Omega_T,
\label{inhomoeqn-des2c}
\\
&w_{\eps m}(\cdot,0) = w^0_\eps \qquad\mbox{ in }V',
\label{inhomoeqn-des2d}
\end{eqnarray}
which is equivalent to the approximate weak equation
\be \label{approx}
\ba{r}
\ds \int\!\!\!\int_{\Omega_T} [(w^0_\eps - w_{\eps m}) D_tv
+ (\bar{\vec z}_{\eps m} + \bar{\vec h}_{\eps m}) \cdot \nabla v] \, dxdt = 0
\\ [4mm]
\forall v\in H^1(0,T;V), v(\cdot, T)=0.
\ea
\ee

By mimicking the procedure of Proposition~\ref{equiv}, it is promptly checked that
\eqref{inhomoeqn-des2b} and \eqref{inhomoeqn-des2c}
may be replaced by the two inequalities
\begin{eqnarray}
&\ds \int\!\!\!\int_{\Omega_T}
[\varphi_\eps(\bar{u}_{\eps m},x) + \varphi_\eps^*(\bar{w}_{\eps m},x)
- \bar{w}_{\eps m} \bar{u}_{\eps m}] \, dxdt\le 0,
\label{mini-inequality1}
\\
&\ds \int\!\!\!\int_{\Omega_T}
[f_{\vec\gamma_\eps}(\nabla \bar{u}_{\eps m}, \bar{\vec z}_{\eps m},x)
- \nabla \bar{u}_{\eps m} \cdot \bar{\vec z}_{\eps m}] \, dxdt \le 0.
\label{mini-inequality2}
\end{eqnarray}

Defining the space
\be
X_{\eps m}= \{(\bar{u}_m, w_m, \bar{\vec z}_m)\in
L^2(0,T; V) \times L^2(\Omega_T)\times L^2(\Omega_T)^N
\mbox{ that fulfill } \eqref{approx}\},
\ee
we conclude that Problem~\ref{problem3} is equivalent to
the following null-minimization problem:

\begin{problem}\label{problem4}
Find $(u_{\eps m}, w_{\eps m}, \vec z_{\eps m}) \in X_{\eps m}$ such that
\be\label{minimization2}
\Phi_\eps(u_{\eps m}, w_{\eps m}, \vec z_{\eps m}) = \inf_{X_{\eps m}} \Phi_\eps= 0.
\ee
\end{problem}

\subsection{A priori estimates}

By the Fenchel inequality \eqref{fenchel}, the inequality \eqref{mini-inequality1}
is tantamount to \eqref{inhomoeqn-des2b}.
By \eqref{approx} and \eqref{inhomoeqn-des2b}
\be
\ba{ll}\label{mini}
&\ds -\int\!\!\!\int_{\Omega_T}
\nabla \bar{u}_{\eps m} \cdot \bar{\vec z}_{\eps m} \, dxdt
= \int_0^T\langle D_tw_{\eps m}, \bar{u}_{\eps m}\rangle_{V',V} \, dt
+ \int\!\!\!\int_{\Omega_T} \nabla \bar{u}_{\eps m} \cdot \bar{\vec h}_{\eps m} \, dxdt
\\
&\ds = \int\!\!\!\int_{\Omega_T} D_t \varphi_\eps^*(w_{\eps m}(x,t),x) \, dxdt
+ \int\!\!\!\int_{\Omega_T} \nabla \bar{u}_{\eps m} \cdot \bar{\vec h}_{\eps m} \, dxdt
\\
&\ds = \int_\Omega
[\varphi_\eps^*(w_{\eps m}(x,T),x) - \varphi_\eps^*(w^0_\eps,x)] \, dx
+ \int\!\!\!\int_{\Omega_T} \nabla \bar{u}_{\eps m} \cdot \bar{\vec h}_{\eps m} \, dxdt.
\ea
\ee
By \eqref{gstar-growth}
\be
\ds \int_\Omega g^*\left(\bar{w}_{\eps m}(x,T),x/\eps\right) \, dx
\ge L \ds\int_\Omega |\bar{w}_{\eps m}(\cdot,T)|^2 -M|\Omega|.
\ee
On the other hand, as the function $f_{\vec\gamma_\eps}$ {\it represents\/} the operator
$\vec\gamma_\eps$ (in the sense of the theory of Fitzpatrick), \eqref{bound-below} yields
\be
\ba{ll}
\ds \int\!\!\!\int_{\Omega_T}
f_{\vec\gamma_\eps}(\nabla \bar{u}_\eps , \bar{\vec z}_\eps,x) \, dxdt
&\ds\ge \int_0^T (\nabla \bar{u}_{\eps m}, \bar{\vec z}_{\eps m}) \, dt
\\ [2mm]
&\ds \geq a\left(\|\nabla \bar{u}_{\eps m}\|^2_{L^2(\Omega)^N} +
\| \bar{\vec z}_{\eps m}\|^2_{L^2(\Omega)^N}\right) -b|\Omega|.
\ea
\ee
By \eqref{mini-inequality2}, then
\be
\ba{ll}
&\ds a\left(\|\nabla \bar{u}_{\eps m}\|^2_{L^2(\Omega)} +
\| \bar{\vec z}_{\eps m}\|^2_{L^2(\Omega)^N}\right) -b|\Omega| + L
\ds\int_\Omega |\bar{w}_{\eps m}(\cdot,T)|^2 -M|\Omega|
\\ [3mm]
&\ds \le \int_\Omega \varphi_\eps^*(w^0_\eps) \, dx -
\int\!\!\!\int_{\Omega_T} \bar{\vec h}_{\eps m} \cdot\nabla \bar{u}_{\eps m} \, dxdt
\\ [3mm]
&\ds \le \int_\Omega \varphi_\eps^*(w^0_\eps) \, dx +
\|\bar{\vec h}_{\eps m}\|_{L^2(\Omega_T)^N} \|\bar{u}_{\eps m}\|_{L^2(0,T;V)}.
\ea
\ee

As in these inequalities one may replace $T$ by any $t\in{} ]0,T]$,
we get the uniform estimates
\be \label{con-est1}
\|\bar{u}_{\eps m}\|_{L^2(0,T;V)} \leq C_1,
\quad \| \bar{\vec z}_{\eps m}\|_{L^2(\Omega_T)^N} \leq C_2,
\ee
where $C_1,C_2,...$ are constants independent of $\eps$.
By the above computation, we also infer that
\be \label{con-est2}
\|\bar{w}_{\eps m}\|_{L^\infty(0,T; H)} \leq C_3,
\ee
and by comparing the terms of \eqref{inhomoeqn-des2a} we
conclude that $w_{\eps m} \in H^1(0,T; V')$ and
\be \label{con-est3}
\|w_{\eps m}\|_{ H^1(0,T; V')} \leq C_4.
\ee

Analogous estimates to \eqref{con-est1} and \eqref{con-est2}
hold for the piecewise interpolate functions, that is,
$w_{\eps m}, u_{\eps m}, {\vec z}_{\eps m}$.
On the other hand, obviously \eqref{con-est3} does not apply to $\bar w_{\eps m}$.

\subsection{Passage to the limit}

On the basis of the above a priori estimates, there exists
$u_\eps, w_\eps,\vec z_\eps$ such that, up to extracting subsequences,
\footnote{ With standard notation, we shall denote the (single-scale)
strong and weak convergence by $\to$ and $\wto$, respectively.
}
\begin{eqnarray}
&u_{\eps m} \wto u_{\eps}
\qquad \mbox{ in }L^2(0,T;V),
\label{con-des1}
\\
&w_{\eps m} \wsto w_{\eps}
\qquad\mbox{ in }L^\infty(0,T;H)\cap H^1(0,T; V'),
\label{con-des2}
\\
&\vec z_{\eps m} \wto \vec z_{\eps}
\qquad\mbox{ in }L^2(\Omega)^N.
\label{con-des3}
\end{eqnarray}
Moreover,
\begin{eqnarray}
&\vec h_{\eps m} \wto \vec h_{\eps}
\qquad \mbox{ in } L^2(\Omega_T)^N,
\label{con-des4}
\\
&w^0_{\eps m} \wto w^0_{\eps}
\qquad\mbox{ in }L^2(\Omega).
\label{con-des5}
\end{eqnarray}

By passing to the limit in \eqref{approx}, we get the equation \eqref{inhomoeqn2a};
namely, $(u_\eps, w_\eps, \vec z_\eps)\in X_\eps$.
Let us next derive \eqref{minimization2} by passing to the inferior limit in
\eqref{minimization1}.
By the weak sequential weak lower semicontinuity of $\varphi_\eps$ and
by \eqref{mini}, we have
\be\label{mini=}
\ba{ll}
&\ds \liminf_{\eps\to 0} -\int\!\!\!\int_{\Omega_T}
\nabla \bar{u}_{\eps m} \cdot \bar{\vec z}_{\eps m} \, dxdt
\\
&\ds \ge \int_\Omega
[\varphi_\eps^*(w_\eps(x,T),x) - \varphi_\eps^*(w^0_\eps,x)] \, dx
+ \int\!\!\!\int_{\Omega_T} \nabla u_\eps \cdot \vec h_\eps \, dxdt
\\
&\ds = \int_0^T \langle D_tw_\eps, u_\eps \rangle_{V',V} \, dt
+ \int\!\!\!\int_{\Omega_T} \nabla u_\eps\cdot \vec h_\eps \, dxdt
\\
&\ds \overset{\eqref{inhomoeqn2a}}{=}
-\int\!\!\!\int_{\Omega_T} \nabla u_\eps \cdot \vec z_\eps \, dxdt.
\ea
\ee
By the weak sequential weak lower semicontinuity of $\varphi_\eps$,
$\varphi_\eps^*$ and $f_\eps$, we then infer that
\begin{eqnarray}
&\ds \int\!\!\!\int_{\Omega_T}
[\varphi_\eps(u_\eps,x) + \varphi_\eps^*(w_\eps,x) - w_\eps u_\eps] \, dxdt\le 0,
\label{mini-inequality1'}
\\
&\ds \int\!\!\!\int_{\Omega_T}
[f_{\vec\gamma_\eps}(\nabla u_\eps, \vec z_\eps,x)
- \nabla u \cdot \vec z_\eps] \, dxdt \le 0,
\label{mini-inequality2'}
\end{eqnarray}
namely
\be
\Phi_\eps(u_\eps, w_\eps, \vec z_\eps) \leq 0;
\ee
that is, $(u_\eps, w_\eps, \vec z_\eps)$ solves Problem~\eqref{problem2}.
We have thus proved the following assertion.

\begin{theorem}\label{existeps}
Let the conditions \eqref{defg}--\eqref{defw} be fulfilled for any fixed $\eps > 0$,
as well as \eqref{bdd}.
The solutions $(u_{\eps m}, w_{\eps m}, \vec z_{\eps m})$ of Problem~\ref{problem3}
then satisfy the uniform estimates \eqref{con-est1}--\eqref{con-est3}.
Therefore there exists $(u_{\eps}, w_{\eps}, \vec z_{\eps})$ such that,
up to extracting subsequences, \eqref{con-des1}--\eqref{con-des3} hold.
\hfill\break\indent
The triplet $(u_{\eps}, w_{\eps}, \vec z_{\eps})$ is then a solution of
Problem~\eqref{problem2} (equivalently, of Problem~\eqref{problem1}).
Finally, the following uniform estimates hold:
\be \label{con-est}
\|u_{\eps}\|_{L^2(0,T;V)},
\| \vec z_{\eps}\|_{L^2(\Omega_T)^N},
\|w_{\eps}\|_{L^\infty(0,T; H)\cap H^1(0,T; V')} \leq \text{Constant.}
\ee
\end{theorem}

\bigskip
\section{Two-scale Formulation}

In this section we introduce two mutually equivalent two-scale formulations,
that we then derive by passing to the limit as $\eps\to 0$ in Problem~\ref{problem1}
(or \ref{problem2}).

We shall denote by $H^1_\sharp(Y)$ the subspace of the functions of $H^1(Y)$ that
have equal traces on opposite faces of $Y$;
these coincide with the restrictions of the $Y$-periodic functions of $H^1(\R^N)$.

We introduce two equivalent two-scale formulation,
in which the constitutive relations are respectively expressed either as inclusions
or as null-minimization principles.

\begin{problem}\label{problem5}
Find
\be\label{2sreg}
\ba{ll}
&u\in L^2(0,T;V), \quad
u_1\in L^2(\Omega_T; H^1_\sharp(Y)),
\\ [3mm]
&w\in L^2(\Omega_T\times Y)\cap H^1(0,T; L^2(Y;V')), \quad
\vec z\in L^2(\Omega_T\times Y)^N,
\ea
\ee
such that
\begin{eqnarray}
&\ba{ll}
&\ds \int\!\!\!\int\!\!\!\int_{\Omega_T\times Y} \big[(w_0 -w) D_tv
+ (\vec z + \vec h)\cdot (\nabla v + \nabla_y v_1) \big] \, dxdtdy =0
\\ [5mm]
&\forall v\in H^1(0,T;V),\, v|_{t=T}=0,
\forall v_1\in L^2(\Omega_T; H^1_\sharp(Y)),
\ea
\label{twoscaleFPsyst2}
\\ [3mm]
&w\in \partial \varphi(u,y) \qquad\mbox{ a.e.\ in }\Omega_T\times Y,
\label{twoscaleFPsyst3}
\\
&\vec z\in \vec\gamma(\nabla u+\nabla_y u_1,y)
\qquad\mbox{ a.e.\ in }\Omega_T\times Y.
\label{twoscaleFPsyst4}
\end{eqnarray}
\end{problem}

Let us next define the space
\be
X_0= \{(u, u_1, w, \vec z) \mbox{ as in $\eqref{2sreg}$ that fulfill }
\eqref{twoscaleFPsyst2}\},
\ee
and the functional
\be \label{Phi-defn}\ba{ll}
&\ds \Phi_0 (u,u_1,w,\vec z) :=
\int\!\!\!\int\!\!\!\int_{\Omega_T\times Y} [\varphi(u,y) + \varphi^*(w,y) -wu
\\ [5mm]
&\ds +f_{\vec\gamma}(\nabla u +\nabla _y u_1, \vec z,y)
- (\nabla u +\nabla _y u_1) \cdot \vec z\,] \, dxdtdy
\qquad\forall (u,u_1,w,\vec z)\in X_0.
\ea\ee 

\medskip
We are now able to introduce our second two-scale formulation.

\begin{problem}\label{problem6}
Find $(u,u_1,w,\vec z)\in X_0$ such that
\be\label{minimization3}
\Phi_0(u,u_1,w,\vec z) = \inf_{X_0} \Phi_0 = 0.
\ee
\end{problem}

\begin{proposition}\label{2s-equiv}
The two-scale Problems~\ref{problem5} and \ref{problem6} are mutually equivalent.
\end{proposition}

\proof
This argument mimics that of Proposition~\ref{equiv}.
The null-minimization of $\Phi_0$ is equivalent to the system of the two inequalities
\begin{eqnarray}
&\ds \int\!\!\!\int\!\!\!\int_{\Omega_T\times Y}
[\varphi(u,y) + \varphi^*(w,y) -wu] \, dxdtdy \le 0,
\label{mini-inequality1=}
\\
&\ds \int\!\!\!\int\!\!\!\int_{\Omega_T\times Y}
[f_{\vec\gamma}(\nabla u +\nabla _y u_1, \vec z,y)
- (\nabla u + \nabla_y u_1) \cdot \vec z\,] \, dxdtdy \le 0,
\label{mini-inequality2=}
\end{eqnarray}
which are respectively equivalent to \eqref{twoscaleFPsyst3} and
\eqref{twoscaleFPsyst4}.
\qed

\medskip
\begin{theorem}\label{2sconv}
Let the assumptions \eqref{defg}--\eqref{defw}, \eqref{bdd} be fulfilled. 
For any $\eps>0$, let $(u_\eps, w_\eps, \vec z_\eps)$ be a solution of
Problem~\ref{problem1} or equivalently of Problem~\ref{problem2}
(this exists by Theorem~\ref{existeps}). Then there exist $u,
w,\vec z$ as in \eqref{2sreg} such that, as $\eps\to 0$ along a
suitable sequence,
\begin{eqnarray}
&u_\eps \w2to u \qquad\mbox{ in }L^2(0,T;V),
\label{2sconv.1}
\\
&\nabla u_\eps \w2to \nabla u + \nabla_y u_1
\quad\text{ in }L^2(\Omega_T \!\times\! Y)^N,
\label{2sconv.2}
\\
&w_\eps \w2to w \qquad\mbox{ in }L^2(\Omega_T \!\times\! Y)^N,
\label{2sconv.3}
\\
&\vec z_\eps \w2to \vec z
\qquad\mbox{ \ in }L^2(\Omega_T\!\times\! Y)^N.
\label{2sconv.4}
\end{eqnarray}
\indent
Moreover, $(u, u_1, w, \vec z)$ is then a solution of Problem~\ref{problem5},
or equivalently of Problem~\ref{problem6}.
\end{theorem}

\proof (i) By Theorem~\ref{existeps} the family of solutions
$\{(u_\eps, w_\eps, \vec z_\eps)\}$ fulfills the uniform estimates
\eqref{con-est}. By Theorems~\ref{cpt-thm} and \ref{twos1} in the
Appendix, then there exist $u, w,\vec z$ as in \eqref{2sreg} that
fulfill \eqref{2sconv.1}--\eqref{2sconv.4} as $\eps\to 0$ along a
suitable sequence. By \eqref{index4} and \eqref{index5}
\begin{eqnarray}
&\vec h_ \eps\2to \vec h
\qquad \mbox{ in } L^2(\Omega_T\!\times\! Y)^N,
\label{2sconv.5}
\\
&w^0_ \eps\2to w^0
\qquad\mbox{ in }L^2(\Omega\!\times\! Y).
\label{2sconv.6}
\end{eqnarray}
By passing to the limit in \eqref{inhomoeqn2a} we then get
the equation \eqref{twoscaleFPsyst2}. 

\medskip
(ii) Next we prove \eqref{twoscaleFPsyst3}.
The null-minimization \eqref{minimization2} is tantamount to
\begin{eqnarray}
&\ds \int\!\!\!\int_{\Omega_T}
[\varphi_\eps(u_\eps,x) + \varphi_\eps^*(w_\eps,x) -w_\eps u_\eps] \, dxdt =0,
\label{in.1}
\\ [2mm]
&\ds \int\!\!\!\int_{\Omega_T} [f_{\vec\gamma_\eps}(\nabla u_\eps, \vec z_\eps,x)
-\nabla u_\eps\cdot \vec z_\eps ] \, dxdt =0.
\label{in.2}
\end{eqnarray}

By \eqref{2sconv.1} and \eqref{2sconv.3},
recalling that $\{w_{\eps}\}$ is also uniformly bounded in $H^1(0,T; V')$, we have
\be
\int\!\!\!\int_{\Omega_T} w_\eps u_\eps \, dxdt \to
\int\!\!\!\int\!\!\!\int_{\Omega_T\times Y} wu \, dxdtdy.
\ee
By \eqref{in.1} and \eqref{2spot.2}, we then infer that
\be\label{in.1=}
\int\!\!\!\int\!\!\!\int_{\Omega_T\times Y}
[\varphi(u,y) + \varphi^*(w,y) -wu] \, dxdtdy \le 0,
\ee
and this is equivalent to \eqref{twoscaleFPsyst3}.

\medskip
(ii) We are left with the proof of \eqref{twoscaleFPsyst4}. 
By \eqref{2spot.2}
\be\label{in.1=}
\liminf_{\eps\to 0} \int\!\!\!\int_{\Omega_T}
f_{\vec\gamma_\eps}(\nabla u_\eps, \vec z_\eps,x) \, dxdtdy
\ge \int\!\!\!\int\!\!\!\int_{\Omega_T\times Y}
f_{\vec\gamma}(\nabla u + \nabla_y u_1, \vec z,y) \, dxdtdy.
\ee

On the other hand, using \eqref{inhomoeqn1a} and \eqref{inhomoeqn2b}
and mimicking \eqref{mini=}, we have
\be
\ba{ll}\label{eq.1}
&\ds -\int\!\!\!\int_{\Omega_T}
\nabla u_\eps \cdot \vec z_\eps \, dxdt
= \int_0^T\langle D_tw_\eps, u_\eps \rangle_{V',V} \, dt
+ \int\!\!\!\int_{\Omega_T} \nabla u_\eps \cdot \vec h_\eps \, dxdt
\\
&\ds = \int\!\!\!\int_{\Omega_T} D_t \varphi_\eps^*(w_\eps(x,t),x) \, dxdt
+ \int\!\!\!\int_{\Omega_T} \nabla u_\eps \cdot \vec h_\eps \, dxdt
\\
&\ds = \int_\Omega
[\varphi_\eps^*(w_\eps(x,T),x) - \varphi_\eps^*(w^0_\eps,x)] \, dx
+ \int\!\!\!\int_{\Omega_T} \nabla u_\eps \cdot \vec h_\eps \, dxdt.
\ea
\ee
By \eqref{2sconv.2}, \eqref{2sconv.3} and \eqref{2spot.2},
\be
\ba{ll}\label{eq.2}
&\ds \liminf_{\eps\to 0} \Big\{ \int_\Omega [\varphi_\eps^*(w_\eps(x,T),x)
- \varphi_\eps^*(w^0_\eps,x)] \, dx
+ \int\!\!\!\int_{\Omega_T} \nabla u_\eps \cdot \vec h_\eps \, dxdt \Big\}
\\
&\ds \ge \int\!\!\!\int_{\Omega\times Y}
[\varphi^*(w(x,T)) - \varphi^*(w^0)] \, dxdy
+ \!\int\!\!\!\int\!\!\!\int_{\Omega_T\times Y}
(\nabla u + \nabla_y u_1) \cdot \vec h \, dxdtdy.
\ea
\ee
Here also we may drop the term in $\nabla_y u_1$.
Moreover, by \eqref{twoscaleFPsyst2} and \eqref{twoscaleFPsyst3},
recalling that $\nabla u$ is independent of $y$,
\be
\ba{ll}\label{eq.3}
&\ds - \int\!\!\!\int\!\!\!\int_{\Omega_T\times Y}
(\nabla u + \nabla_y u_1) \cdot (\vec z +\vec h) \, dxdtdy 
\\[4mm]
&\ds = \int_Y \! dy \!\! \int_0^T \! \langle D_tw, u \rangle_{V',V} \, dt
= \int\!\!\!\int\!\!\!\int_{\Omega_T\times Y} D_t \varphi^*(w(x,y,t)) \, dxdtdy
\\[4mm]
&\ds = \int\!\!\!\int_{\Omega\times Y} [\varphi^*(w(x,y,T)) - \varphi^*(w^0(x,y))] \, dxdy.
\ea
\ee
By \eqref{eq.1}, using \eqref{eq.2} and \eqref{eq.3}, we have
\be
\liminf_{\eps\to 0} \; -\int\!\!\!\int_{\Omega_T}
\nabla u_\eps\cdot \vec z_\eps \, dxdt
\ge - \int\!\!\!\int_{\Omega_T} \nabla u \cdot \vec z \, dxdt.
\ee

By passing to the inferior limit in \eqref{in.2} and using \eqref{2spot.2},
we then get
\be\label{in.2=}
\int\!\!\!\int\!\!\!\int_{\Omega_T\times Y}
[f_{\vec\gamma}(\nabla u + \nabla_y u_1, \vec z,y)
- (\nabla u + \nabla_y u_1)  \cdot \vec z\,] \, dxdtdy \le 0,
\ee
which is tantamount to \eqref{twoscaleFPsyst4}.
\qed

\bigskip
\section{Single-Scale Formulation (Homogenization)} 

In this section we derive a single-scale formulation
(i.e., a homogenized problem) from the two equivalent
two-scale Problems~\ref{problem5} and \ref{problem6},
and prove a homogenization theorem.
Along the lines of the previous sections, we introduce two equivalent formulations,
in which the constitutive relations are respectively expressed either as inclusions
or as null-minimization principles.

Let the convex function $\varphi_0$ and the maximal monotone map $\vec\gamma_0$
be respectively defined as in Propositions~\ref{upscc} and \ref{upscmm}.
Here is our first single-scale formulation.

\begin{problem}\label{problem7}
Find
\be\label{2sreg'}
u\in L^2(0,T;V), \quad
w\in L^2(\Omega_T)\cap H^1(0,T;V'), \quad
\vec z\in L^2(\Omega_T)^N,
\ee
such that
\begin{eqnarray}
&\ba{rr}
&\ds \int\!\!\!\int_{\Omega_T} \big[(w_0 -w) D_tv
+ (\vec z + \vec h)\cdot \nabla v \big] \, dxdt =0
\\ [5mm]
&\forall v\in H^1(0,T;V),\, v|_{t=T}=0,
\ea
\label{twoscaleFPsyst2'}
\\ [3mm]
&w\in \partial \varphi_0(u) \qquad\mbox{ a.e.\ in }\Omega_T,
\label{twoscaleFPsyst3'}
\\
&\vec z\in \vec\gamma_0(\nabla u)
\qquad\mbox{ a.e.\ in }\Omega_T.
\label{twoscaleFPsyst4'}
\end{eqnarray}
\end{problem}

We already know that the weak equation \eqref{twoscaleFPsyst2'} is equivalent to
the PDE
\be \label{PDE+}
D_tw - \nabla\cdot \vec z = \nabla\cdot \vec h
\qquad\mbox{ in $V'$, a.e.\ in }]0,T[,
\ee
coupled with the initial condition
\be
w(\cdot,0)= w^0 \qquad\mbox{ a.e.\ in }\Omega.
\ee

Let us next define the space
\be
X_0= \{(u, w, \vec z) \mbox{ as in $\eqref{2sreg'}$ that fulfill }
\eqref{twoscaleFPsyst2'}\},
\ee
the mutually orthogonal spaces
\be\label{ortho.1}
\ba{ll}
&W = \{\nabla\phi: \phi\in W^{1,p}(Y)\},
\\[2mm]
&Z = \{\vec v\in L^{p'}_*(Y): \int_Y w(y) \, dy=0, \nabla \cdot \vec v=0\},
\ea
\ee
and the functionals
\be\label{ort}
\ba{rr}
\ds F_0(\vec \xi,\vec \eta) = \inf_{\vec v\in W,\vec w\in Z}
\int_Y f_{\vec\gamma}(\vec\xi + \vec v(y),\vec\eta + \vec w(y),y)  \, dy
\qquad \forall \vec \xi,\vec \eta\in \R^N,
\ea
\ee
\be \label{Phi-defn}
\ba{rr}
\ds \Phi_0 (u, w,\vec z) :=
\int\!\!\!\int_{\Omega_T}
[\varphi_0(u) + \varphi_0^*(w) -wu + F_0(\nabla u, \vec z)
- \nabla u \cdot \vec z\,] \, dxdt
\\[3mm]
\forall (u, w,\vec z)\in X_0.
\ea\ee

We are now able to introduce our second single-scale formulation.

\begin{problem}\label{problem8}
Find $(u, w,\vec z)\in X_0$ such that
\be\label{minimization3'}
\Phi_0(u,u_1,w,\vec z) = \inf_{X_0} \Phi_0 = 0.
\ee
\end{problem}

\begin{proposition}\label{2s-equiv'}
The single-scale Problems~\ref{problem7} and \ref{problem8} are mutually equivalent.
\end{proposition}

\proof
This argument mimics that of Proposition~\ref{equiv}.
The null-minimization of $\Phi_0$ is equivalent to the system of the two inequalities
\begin{eqnarray}
&\ds \int\!\!\!\int_{\Omega_T} [\varphi_0(u) + \varphi_0^*(w) -wu] \, dxdt \le 0,
\label{mini-inequality1+}
\\
&\ds \int\!\!\!\int_{\Omega_T}
[F_0(\nabla u, \vec z) - \nabla u \cdot \vec z\,] \, dxdt \le 0,
\label{mini-inequality2+}
\end{eqnarray}
which are respectively equivalent to \eqref{twoscaleFPsyst3'} and
\eqref{twoscaleFPsyst4'}.
\qed

\medskip
We shall use the two-scale decomposition
\be
\ba{ll}
&\ds \widehat{u}(x) := \int_Y u(x,y) \, dy
\\ [4mm]
&\widetilde{u}(x,y) := u(x,y) - \widehat{u}(x)
\ea
\qquad\mbox{ for a.e.\ }(x,y) \in \Omega\times Y.
\ee

\begin{proposition}\label{upscaling}
If $(u,u_1,w,\vec z)$ is a solution of Problem~\ref{problem5}
or equivalently of Problem~\ref{problem6}
(such a solution exists by Theorem~\ref{2sconv}), then
$(u,\widehat{w},\widehat{\vec z})$ is a solution of Problem~\ref{problem7}
or equivalently of Problem~\ref{problem8}.
\end{proposition}

\proof
Selecting either $v=0$ or $v_1=0$ in the equation \eqref{twoscaleFPsyst2},
we respectively get
\be
\ba{rr}
&\ds \int\!\!\!\int_{\Omega_T} \big[(\widehat{w_0} - \widehat{w}) D_tv
+ (\widehat{\vec z} + \widehat{\vec h}) \cdot \nabla v \big] \, dxdt =0
\\ [4mm]
&\forall v\in H^1(0,T;V),\, v|_{t=T}=0,
\ea
\ee
\be
\int\!\!\!\int\!\!\!\int_{\Omega_T\times Y}
(\widetilde{\vec z} + \widetilde{\vec h})\cdot \nabla_y v_1 \, dxdtdy =0
\qquad\forall v_1\in L^2(\Omega_T; H^1_\sharp(Y)).
\ee
These integral equations respectively correspond to the following
coarse- and fine-scale PDEs:
\begin{eqnarray}
&D_t\widehat{w} - \nabla \cdot \widehat{\vec z} = \nabla \cdot \widehat{\vec h}
\qquad\mbox{ in $V'$, a.e.\ in }]0,T[,
\\
&- \nabla_y \cdot \widetilde{\vec z} = \nabla_y \cdot \widetilde{\vec h}
\qquad\mbox{ in $H^1_\sharp(Y)'$, a.e.\ in }\Omega_T.
\end{eqnarray}

By Propositions~\ref{upscc} and \ref{upscmm},
the single-scale constitutive relations
\eqref{twoscaleFPsyst3'} and \eqref{twoscaleFPsyst4'} follow from
\eqref{twoscaleFPsyst3} and \eqref{twoscaleFPsyst4}.
\qed

\medskip
\begin{theorem}\label{1sconv}
Let the assumption \eqref{defg}--\eqref{defw}, \eqref{bdd} be fulfilled.
For any $\eps>0$, let $(u_\eps, w_\eps, \vec z_\eps)$ be a solution of
Problem~\ref{problem1} or equivalently of Problem~\ref{problem2}
(this exists by Theorem~\ref{existeps}). Then there exist $u,
w,\vec z$ as in \eqref{2sreg'} such that, as $\eps\to 0$ along a
suitable sequence,
\begin{eqnarray}
&u_\eps \wto u \qquad\mbox{ in }L^2(0,T;V),
\label{2sconv.1'}
\\
&w_\eps \wto w \qquad\mbox{ in }L^\infty(0,T; H)\cap H^1(0,T;V'),
\label{2sconv.3'}
\\
&\vec z_\eps \wto \vec z
\qquad\mbox{ \ in }L^2(\Omega_T)^N.
\label{2sconv.4'}
\end{eqnarray}
\indent
This entails that $(u, w, \vec z)$ is a solution of the homogenized
Problem~\ref{problem7}, or equivalently of Problem~\ref{problem8}.
\end{theorem}

\bigskip
\section{Appendix}

Here we briefly review the notion of two-scale convergence, and
some related properties of integral functionals.

\subsection{Two-scale convergence}

This notion was introduced by Nguetseng \cite{Ng},
and was further developed by Allaire and others, see e.g.\ \cite{Al},\cite{LuNgWa}.

Let us denote by $Y= {}]0,1[^N$ the fundamental periodicity-cell, and
by $\eps >0$ a small parameter which we shall let eventually vanish.
Let us fix any $p\in {}]1,+\infty[$ and define the conjugate index $p':= p/(p-1)$.
Let us denote by $C_ \sharp(Y)$ ($W^{1,p}_ \sharp(Y)$, resp.)
the space of continuous ($W^{1,p}_\sharp$, resp.) functions
$\R^N\to \R$ that are $Y$-periodic and have equal traces on opposite faces of $Y$.
By the index $*$ we shall denote subspaces of functions with vanishing average:
e.g., $L^1_*(Y) = \{w\in L^1(Y): \int_Y w(y) \, dy=0 \}$.

\begin{definition}[Weak two-scale convergence]\label{defn-w2scale}
We shall say that a sequence $\{u_\eps\}$ of functions in $L^p(\Omega)$
{\it weakly two-scale converges\/} to a limit function $u\in L^p(\Omega\times Y)$,
and write $u_n\w2to u$,
whenever
\be\label{defn-2scale-eqn1}
\ds{\int_\Omega u_\eps (x)\phi(x,x/\eps) \, dx \to
\int\!\!\!\int_{\Omega\times Y} u(x,y)\phi(x,y) \, dxdy}
\qquad\forall\phi \in L^{p'}(\Omega ;C_ \sharp(Y)).
\ee
\end{definition}

For instance,
\be\label{weaktwo}
x\sin (2\pi x/\eps) \w2to x\sin (2\pi y)
\qquad
\mbox{ in }L^p(\Omega\times Y), \forall p\in {}]1,+\infty[.
\ee
Notice that the weak two-scale limit is unique, if it exists.

This definition is trivially extended to time-dependent functions.
For any $p,r\in {}]1,+\infty[$, we shall say that a family $\{u_\eps\}$ of functions in
$L^r(0,T, L^p(\Omega))$ {\it weakly two-scale converges\/}
to a limit $u\in L^r(0,T, L^p(\Omega\times Y))$ whenever
\be\label{defn-2scale-eqn2}
\begin{array}{r}
\ds{\int\!\!\!\int_{\Omega_T} u_\eps (x,t)\phi(x,x/\eps,t) \, dxdt \to
\int\!\!\!\int\!\!\!\int_{\Omega_T\times Y} u(x,y,t)\phi(x,y,t) \, dxdydt}
\\ [6mm]
\qquad\forall\phi \in L^{r'}(0,T, L^{p'}(\Omega ;C_ \sharp(Y))).
\end{array}
\ee

The results that follow also trivially take over to time-dependent functions.

\begin{theorem}\label{cpt-thm}
If $\{u_\epsilon\}$ is a bounded sequence in $L^p( \Omega)$ ($p\in {}]1,+\infty[$),
then there exists $u\in L^p(\Omega\times Y)$ such that,
as $\eps\to 0$ along a suitable subsequence,
$u_\eps \w2to u$ in $L^p(\Omega\times Y)$.
\end{theorem}

\begin{proposition}
If $u_\eps\to u$ in $L^p(\Omega)$ ($p\in {}]1,+\infty[$),
then $u_\eps \w2to u$ ($=u(x)$) in $L^p(\Omega\times Y)$.
\hfil\break\indent
On the other hand, if $u_\eps \w2to u$ ($=u(x,y)$) in $L^p(\Omega\times Y)$, then
the sequence $\{u_\eps\}$ is bounded in $L^p(\Omega)$, and
$u_\eps \wto \int_Y u(x,y) \, dy$ in $L^p(\Omega )$.
\end{proposition}

For any measurable function $u:\Omega\times Y\to \R$ such that
$u(x,\cdot)\in L^1(Y)$ for a.e.\ $x\in\Omega$,
we define the average component $\widehat{u}$
and the fluctuating component $\widetilde{u}$ as follows:
\be\label{decomp}
\ba{ll}
&\ds \widehat{u}(x) := \int_Y u(x,y) \, dy
\\ [4mm]
&\widetilde{u}(x,y) := u(x,y) - \widehat{u}(x)
\ea
\qquad\mbox{ for a.e.\ }(x,y) \in \Omega\times Y.
\ee
Thus $\widetilde{u}(x,\cdot)\in L^1_*(Y)$ for a.e.\ $x\in\Omega$.

\begin{proposition}
If $\{u_\eps\}$ is a sequence  in $L^p(\Omega)$ ($p\in {}]1,+\infty[$) that
two-scale converges to $u\in L^p(\Omega\times Y)$, then
\begin{equation}\label{2scalelowerlimit} \liminf_{\eps \to 0}
\|u_\eps \|_{L^p(\Omega)}\, \geq\,
  \|u\|_{L^p(\Omega \times Y)} \,\geq\, \|\widehat{u}\|_{L^p(\Omega)}.
\end{equation}
\end{proposition}

\begin{definition}[Strong two-scale convergence]\label{defn-s2scale}
We shall say that a sequence $\{u_\eps\}$ of functions in $L^p(\Omega)$
($p\in {}]1,+\infty[$)
strongly two-scale converges to $u=u(x,y)$ in $L^p(\Omega\times Y)$,
and write $u_\eps \2to u$, whenever
$$
u_\eps \w2to u \quad\text{ in }L^p(\Omega\times Y),
\qquad\mbox{and}\qquad
\|u_\eps \|_{L^p(\Omega)} \to \|u\|_{L^p(\Omega \times Y)}.
$$
\end{definition}

For instance the sequence in \eqref{weaktwo} is strongly two-scale convergent,
whereas
$x\sin (2\pi x/\eps) + x\sin (2\pi x/\eps^2)$ is just weakly two-scale convergent
to $x\sin (2\pi y)$.

The next result is one of the major tools for the application of two-scale convergence
to the homogenization of PDEs.

\begin{theorem}\label{twos1}
Let $\{u_\eps\}$ be a bounded sequence in $W^{1,p}(\Omega)$ ($p\in {}]1,+\infty[$).
Then there exist
$(u,u_1) \in W^{1,p}(\Omega) \times L^p(\Omega, W^{1,p}_\sharp(Y))$
such that, as $\eps\to 0$ along a suitable subsequence,
\begin{equation}\label{twosclimit1}
\left\{\begin{array}{ll}
u_\eps\wto u \quad\text{ in } W^{1,p}(\Omega),
\\ [2mm]
u_\eps \w2to u \quad\text{ in } L^p(\Omega),
\\ [2mm]
\nabla u_\eps \w2to \nabla u + \nabla_y u_1 \quad\text{ in } L^p(\Omega)^N.
\end{array}\right.
\end{equation}
\end{theorem}

\subsection{On the measurability of multivalued mappings}

Let us assume that $(S,\mathcal{A})$ is a measurable space and
that $B$ is a separable and reflexive real Banach space with dual $B'$.
We remind the reader that a multivalued mapping $g: S\to \mathcal{P}(B')$ is
called {\it measurable\/} if
\be\label{meas1}
\ba{ll}
\ds g^{-1}(R) := \big\{x\in S: g(x)\cap R\not= \emptyset \big\}
\\ [2mm]
\mbox{is measurable, for any open set $R\subset B'$. }
\ea
\ee
By a classical theorem of Pettis, see e.g.\ \cite{Yo}, it is equivalent to refer to measurability
with respect to the weak or to the strong topology of the separable space $B'$.

Moreover $g$ is called {\it closed-valued\/} if $g(s)$ is closed for a.e.\ $s\in S$.
It is known that if $g$ is closed-valued and measurable, then it has a measurable selection,
see e.g.\ Sect.~III.2 of \cite{CaVa} or Sect.~8.1 of \cite{IoTi}.
This means that there exists a measurable single-valued mapping $f: S\to B'$ such that
$f(x)\in g(x)$ for a.e.\ $x\in S$.

For any domain $D\subset \R^N$, let us denote by $\mathcal{L}(D)$ and
$\mathcal{B}(D)$ the $\sigma$-algebras of Borel- and Lebesgue-measurable functions
$D\to \R$, respectively.
Let us also denote by ${\cal B}(B) \!\otimes\! {\cal L}(Y)$
the $\sigma$-algebra generated by the Cartesian product
of the Lebesgue and Borel $\sigma$-algebras ${\cal B}(B)$ and ${\cal L}(Y)$.

\begin{lemma}
Let us assume that
\be
\ba{lll}\label{meas2}
&\alpha: B\!\times\! Y \to \mathcal{P}(B')
\mbox{ \ is measurable w.r.t. }\mathcal{B}(B)\!\otimes\! \mathcal{L}(Y),
\\ [2mm]
&\alpha(\xi,y) \mbox{ \ is closed for any $\xi$ and a.e.\ }y.
\ea
\ee
For any $\mathcal{L}(Y)$-measurable mapping $v: Y\to B$,
the multivalued mapping $\beta: y\mapsto\alpha(v(y),y)$
is then closed-valued and measurable.
\end{lemma}

\proof
Let us set $\gamma_v(y) =(v(y),y)$ for any $y\in Y$, so that
$\beta =\alpha \circ\gamma_v$ in $Y$.
The mapping $\gamma_v: Y\to B \!\times\! Y$ is clearly measurable.
Because of \eqref{meas2}$_2$ the set $\beta(y)$ is closed for a.e.\ $y$.
For any open set $R\subset B'$, by \eqref{meas2}$_1$
$\alpha^{-1}(R)\in \mathcal{B}(B)\!\otimes\! \mathcal{L}(Y)$.
By the ${\cal L}(Y)$-measurability of $v$, we conclude that
$\beta ^{-1}(R) = \gamma_v^{-1} (\alpha^{-1}(R))
= \big\{y\in Y: (v(y),y) \in \alpha^{-1}(R)\big\} \in {\cal L}(Y)$.
\qed

\subsection{Two-scale limit of integral functionals}

\begin{proposition}\label{convex-lsc} \cite{ViCalVar07}
(i) If $\phi: \R^N\times Y \to \R\cup \{+\infty\}$ is
$\mathcal{B}(\R^N) \bigotimes\mathcal{L}(Y)$-measurable,
then for any measurable function $v: \Omega\to \R^N$,
the mappings $x\mapsto\phi\left(v(x),x/\eps \right)$ and
$(x,y)\mapsto\phi(v(x,y),y)$ are measurable.
\hfil\break\indent
(ii) Let $\phi$ be also convex with respect to the first variable for a.e.\ $y$,
and assume that there exist $C\in\R^N$ and $h\in L^1(Y)$ such that
\be\label{lower-boundphi.1}
\phi(\vec v,y )\geq C\cdot \vec v + h(y)
\qquad
\forall \vec v\in\R^N, \mbox{ for a.e.\ }y\in Y.
\ee
Let us define the functionals $\Psi_\eps: L^p(\Omega)^N \to \R\cup\{+\infty\}$
and $\Psi: L^p(\Omega\times Y)^N \to \R\cup\{+\infty\}$ by
\begin{eqnarray}
&\ds \Psi_\eps (\vec v) := \int_\Omega\phi (\vec v(x),x/\eps) \, dx
\qquad\forall \vec v\in L^p(\Omega)^N,
\label{lower-boundphi.2}
\\ [2mm]
&\ds \Psi(\vec v)  := \int\!\!\!\int_{\Omega\times Y}\phi (\vec v(x,y),y) \, dxdy
\qquad\forall \vec v\in L^p(\Omega\times Y)^N.
\label{lower-boundphi.3}
\end{eqnarray}
\indent
These functionals are well-defined, convex and lower semicontinuous, respectively
in $L^p(\Omega)^N$ and $L^p(\Omega\times Y)^N$.
\hfil\break\indent
(iii) Under the above assumptions, for any sequence
$\{\vec v_\eps\}$ in $L^p(\Omega)^N$,
\begin{eqnarray}
&\vec u_\eps \2to \vec u \mbox{ \ in }L^p(\Omega\times Y)^N
\quad\Rightarrow\quad
\Psi_\eps(\vec v_\eps) \to \Psi(\vec v),
\label{2spot.1}
\\ [2mm]
&\vec u_\eps \w2to \vec u \mbox{ \ in }L^p(\Omega\times Y)^N
\quad\Rightarrow\quad
\liminf_{\eps\mapsto 0} \Psi_\eps(\vec v_\eps) \geq \Psi(\vec v).
\label{2spot.2}
\end{eqnarray}
\end{proposition}

\medskip
It is known that the convex conjugate functionals $\Psi_\eps^*$ and $\Psi^*$
then coincide with the integral functionals of the convex conjugate
of the respective integrands.

\subsection{Scale-integration of cyclically maximal monotone operators}

Let us first set
\be\label{ortho.1}
W = \{\nabla\phi: \phi\in W^{1,p}(Y)\},
\qquad
Z = \{\vec v\in L^{p'}_*(Y): \nabla \cdot \vec v=0\},
\ee
and notice the following orthogonality relation:
\be\label{orthog}
\int_Y \vec u(y) \cdot \vec w(y) \, dy =0
\qquad\forall \vec u\in W, \forall \vec w\in Z.
\ee

Let us assume that
\be
\ba{ll}
&\varphi: \R^N \!\times\! Y \to \R^N
\mbox{ \ is measurable w.r.t. ${\cal B}(\R^N)\!\otimes\! {\cal L}(Y)$, }
\\[2mm]
&\exists p\in{} ]1,+\infty[, \exists c_1,...,c_4>0:
\\[2mm]
&c_1|\vec \xi|^p - c_2 \le \varphi(\vec \xi,y) \le c_3|\vec \xi|^p + c_4
\qquad \forall \vec \xi\in \R^N,\mbox{ for a.e.\ }y\in Y,
\ea
\ee
and set
\begin{eqnarray}
&\ds \varphi_0(\vec \xi) = \inf_{\vec v\in W}
\int_Y \varphi(\vec \xi+ \vec v(y),y) \, dy
\qquad\forall \vec \xi\in \R^N,
\label{ortho.2}
\\
&\ds \psi_0(\vec \eta) = \inf_{\vec v\in Z}
\int_Y \varphi^*(\vec \eta+ \vec v(y),y) \, dy
\qquad\forall \vec \eta\in \R^N.
\label{ortho.3}
\end{eqnarray}

\begin{proposition}\label{upscc} \cite{ViCalVar09} 
The function $\psi_0$ is the convex conjugate of $\varphi_0$.
\hfill\break\indent
If $\vec u\in L^p(Y)^N$ and $\vec w\in L^{p'}(Y)^N$ are such that
\begin{eqnarray}
&\ds \widetilde{\vec u}\in W, \quad \widetilde{\vec w}\in Z,
\\
&\ds \vec w(y) \in\partial \varphi(\vec u(y),y) \qquad\mbox{ for a.e.\ }y\in Y,
\end{eqnarray}
then
\begin{eqnarray}
&\widehat{\vec w}\in\partial \varphi_0(\widehat{\vec u}),
\qquad
\widehat{\vec u}\in\partial \psi_0(\widehat{\vec w}),
\\
&\ds \varphi_0 (\widehat{\vec u}) = \int_Y \varphi(\vec u(y),y) \, dy,
\qquad
\psi_0(\widetilde{\vec w})
= \int_Y \varphi^*(\vec w(y),y) \, dy.
\end{eqnarray}
\end{proposition}

This result takes over to noncyclically maximal monotone operators.

\subsection{Scale-integration of maximal monotone operators}

Let us assume that
\be
\ba{ll}
&\vec\gamma: \R^N \!\times\! Y \to \mathcal{P}(\R^N)
\mbox{ \ is measurable w.r.t. $\mathcal{B}(\R^N) \otimes \mathcal{L}(Y)$, }
\\[2mm]
&\mbox{is maximal monotone w.r.t.\ the first argument for a.e.\ } y\in Y,
\\[2mm]
&\mbox{and is represented by a function $f_{\vec\gamma}(\cdot,y)$ for a.e.\ } y\in Y.
\ea
\ee
Let us then set
\be\label{ort}
\ba{rr}
\ds F_0(\vec \xi,\vec \eta) = \inf_{\vec v\in W,\vec w\in Z}
\int_Y f_{\vec\gamma}(\vec\xi + \vec v(y),\vec\eta + \vec w(y),y)  \, dy
\qquad \forall \vec \xi,\vec \eta\in \R^N.
\ea
\ee

As $f_{\vec\gamma}$ is a representative map, for any $\vec \xi,\vec \eta,\vec v,\vec w$
as above we have
\be
\ba{ll}\label{ort}
\ds \int_Y f_{\vec\gamma}(\vec\xi + \vec v(y),\vec\eta + \vec w(y),y)  \, dy
\ge \int_Y [\vec\xi + \vec v(y)] \cdot [\vec\eta + \vec w(y)]  \, dy
\overset{\eqref{orthog}}{=} \vec\xi \cdot \vec\eta.
\ea
\ee
By taking the infimum with respect to $\vec v\in W$ and $\vec w\in Z$ we thus get
\be
F_0(\vec \xi,\vec \eta) \ge \vec\xi \cdot \vec\eta
\qquad \forall \vec \xi,\vec \eta\in \R^N.
\ee
so that $F_0$ is also a representative function.

\begin{proposition}\label{upscmm} \cite{ViMon}
The function $F_0$ represents a maximal monotone map
$\vec\gamma_0: \R^N \to \mathcal{P}(\R^N)$. 
If
$\vec u\in L^p(Y)^N$ and $\vec w\in L^{p'}(Y)^N$ are such that
\begin{eqnarray}
&\ds \widetilde{\vec u}\in W, \quad \widetilde{\vec w}\in Z,
\\
&\ds \vec w\in \vec\gamma(\vec u,y) \qquad\mbox{ for a.e.\ }y\in Y,
\end{eqnarray}
then
\be
\widehat{\vec w}\in \vec\gamma_0(\widehat{\vec u}).
\ee
\end{proposition}

\bigskip\bigskip
\centerline{\bf Acknowledgment}
\medskip

This work was initiated during the visit of A.V.\ at the Department of Mathematics the Indian Institute of Science in Bangalore,
and was continued during the visit of A.K.N.\ at the Dipartimento di Matematica
dell'Universit\`a di Trento.
Support and hospitality of both institutions are warmly acknowledged.

The work of A.V.\ was also partially supported by the MIUR-PRIN '10-11 grant for the project
``Calculus of Variations".

\end{document}